\documentclass[12pt]{amsart}
\usepackage{amsmath, amssymb, latexsym}
\usepackage{mathrsfs}
\usepackage{enumerate}
\usepackage{color}
\usepackage{graphics}
\usepackage{verbatim}

\usepackage[all, knot]{xy}
\xyoption{arc}

\newtheorem{Thm}{Theorem}[section]
\newtheorem{Proposition}[Thm]{Proposition}
\newtheorem{Corr}[Thm]{Corollary}
\newtheorem{Lem}[Thm]{Lemma}
\newtheorem{Def}[Thm]{Definition}

\theoremstyle{definition}

\newcommand{\nc}{\newcommand}

\nc{\Lemma}{\begin{Lem}}
\nc{\enlemma}{\end{Lem}}
\nc{\Cor}{\begin{Corr}}
\nc{\encor}{\end{Corr}}
\nc{\Th}{\begin{Thm}}
\nc{\entheorem}{\end{Thm}}
\nc{\Prop}{\begin{Proposition}}
\nc{\enprop}{\end{Proposition}}

\numberwithin{equation}{section}

\renewcommand{\to}[1][]{\xrightarrow{#1}{}}

\newcommand{\Z}{\mathbf{Z}}
\newcommand{\Q}{\mathbf{Q}}
\newcommand{\R}{\mathbf{R}}

\newcommand{\A}{\mathbf{A}}

\newcommand{\g}{\mathfrak{g}}
\newcommand{\gsl}{\mathfrak{sl}}

\nc{\on}{\operatorname}

\newcommand{\Hom}{\on{Hom}}

\newcommand{\End} { {\rm End}}
\newcommand{\id} { {\rm id}}

\nc{\cor}{\mathbf{k}}
\nc{\KLR}{Khovanov-Lauda-Rouquier algebra}
\nc{\KLRs}{Khovanov-Lauda-Rouquier algebras}
\nc{\seteq}{\mathbin{:=}}
\newcommand{\soplus}{\mathop{\mbox{\normalsize$\bigoplus$}}\limits}
\nc{\cl}{\colon}
\nc{\set}[2]{\left\{#1\mid #2\right\}}
\nc{\Id}{\operatorname{Id}}
\nc{\Ker}{\on{Ker}}
\nc{\Coker}{\on{Coker}}
\nc{\coh}{\mathrm{coh}}
\nc{\Mod}{\on{Mod}}
\nc{\Modc}{\on{Mod_\coh}}
\nc{\Proj}{\on{Proj}}
\nc{\Rep}{\on{Rep}}
\newcommand{\isoto}[1][]{\mathop{\xrightarrow[#1]%
{{\raisebox{-.6ex}[0ex][-.6ex]{$\mspace{2mu}\sim\mspace{2mu}$}}}}}

\nc{\To}[1][\quad]{\to[\;#1\;]}
\nc{\hs}{\hspace*}
\nc{\vs}{\vspace*}
\nc{\bF}{\overline{F}}
\nc{\epi}{\twoheadrightarrow}
\nc{\mono}{\rightarrowtail}
\nc{\be}{\begin{enumerate}}
\nc{\ee}{\end{enumerate}}
\nc{\ba}{\begin{array}}
\nc{\ea}{\end{array}}
\nc{\eq}{\begin{eqnarray}}
\nc{\eneq}{\end{eqnarray}}
\nc{\eqn}{\begin{eqnarray*}}
\nc{\eneqn}{\end{eqnarray*}}
\nc{\ran}{\rangle}
\nc{\lan}{\langle}
\nc{\bl}{\bigl(}
\nc{\br}{\bigr)}
\nc{\bnum}{\be[{\rm(i)}]}
\nc{\enum}{\ee}
\nc{\bna}{\be[{\rm(a)}]}
\nc{\Proof}{\begin{proof}}
\nc{\QED}{\end{proof}}
\newcommand{\scbul}{{\,\raise1pt\hbox{$\scriptscriptstyle\bullet$}\,}}
\nc{\tens}{\mathop\otimes}
\nc{\E}[1][i]{{\mathsf{E}_{#1}^\Lambda}}
\nc{\F}[1][i]{{\mathsf{F}_{#1}^\Lambda}}
\nc{\x}[1][{1}]{a^\Lambda(x_{#1})}
\nc{\vphi}{\varphi}
\nc{\haut}{\mathrm{ht}}
\nc{\al}{\alpha}
\nc{\La}{\Lambda}
\nc{\Gr}{\on{Gr}}
\nc{\la}{\lambda}
\nc{\noi}{\noindent}
\nc{\eps}{\varepsilon}
\nc{\RL}{R^\La}
\nc{\One}{\mathbf{1}}
\nc{\comout}{}
\nc{\bigwr}{\mbox{\large$\wr$}}
\nc{\heps}{\widehat{\eps}}
\nc{\heta}{\widehat{\eta}}
\nc{\tE}{\widetilde{E}}
\nc{\Fil}{\Gamma}
\nc{\op}{\mathrm{opp}}
\nc{\shc}{\mathscr{C}}
\nc{\Fct}{\on{Fct}}
\nc{\HH}{\mathsf{H}}
\nc{\nn}{\nonumber}
\renewcommand{\Im}{\on{Im}}
\newlength{\my}
\begin{document}

\title[Biadjointness in cyclic Khovanov-Lauda-Rouquier Algebras]
{Biadjointness in cyclic Khovanov-Lauda-Rouquier Algebras}
\author[Masaki Kashiwara]{Masaki Kashiwara}
\thanks{This work was supported by Grant-in-Aid for
Scientific Research (B) 22340005, Japan Society for the Promotion of Science.}
\address{Research Institute for Mathematical Sciences, Kyoto University, Kyoto 606-8502, Japan}
\email{masaki@kurims.kyoto-u.ac.jp}

\date{\today}


\begin{abstract}
In this paper, we prove that a pair of functors $\E$ and $\F$ appearing
in the categorification of irreducible highest weight modules of
quantum groups via cyclotomic Khovanov-Lauda-Rouquier algebras
is a biadjoint pair.
\end{abstract}

\maketitle


\section{Introduction}

Lascoux-Leclerc-Thibon (\cite{LLT}) conjectured that
the irreducible representations of Hecke algebras of type $A$ are
controlled by the upper global basis (\cite{Kash91,Kash93})
(or dual canonical basis (\cite{Lus93}) of
the basic representation of the affine quantum group $U_q(A^{(1)}_\ell)$.
Then Ariki (\cite{A}) proved this conjecture by generalizing it to 
cyclotomic affine Hecke algebras.
The crucial ingredient there 
was the fact that the cyclotomic affine Hecke algebras
categorify the irreducible highest weight representations
of $U(A^{(1)}_\ell)$. Because of the lack of grading on
the cyclotomic affine Hecke algebras, these algebras do not
categorify the representation of the quantum group.

Then Khovanov-Lauda and Rouquier 
introduced independently
a new family of graded algebras, a generalization of
affine Hecke algebras of type $A$,
in order to categorify arbitrary quantum groups (\cite{KL09, KL08, R08}).
These algebras are  called  {\em \KLRs} or 
{\em quiver Hecke algebras.}

 Let $U_q(\g)$ be the quantum
group associated with a symmetrizable Cartan datum and let
$\{R(\beta)\}_{\beta \in Q^{+}} $ be the corresponding
Khovanov-Lauda-Rouquier algebras. Then it was shown in \cite{KL09,
KL08} that there exists an algebra isomorphism
$$U_{\A}^{-}(\g) \simeq \bigoplus_{\beta \in
Q^{+}}K\bl\Proj(R(\beta))\br,$$ where $U_{\A}^-(\g)$ is the integral form
of the half $U_q^-(\g)$ of $U_{q}(\g)$ with $\A = \Z[q, q^{-1}]$,
and $K\bl\Proj(R(\beta))\br$ is the Grothendieck group of finitely generated
projective graded $R(\beta)$-modules. Moreover, when the generalized Cartan
matrix is a symmetric matrix, Varagnolo and Vasserot proved that
{\em lower global basis} introduced by the author or Lusztig's {\it canonical
basis} corresponds to the isomorphism classes of indecomposable
projective $R$-modules  under this isomorphism (\cite{VV09}).

For each dominant integral weight $\Lambda \in P^{+}$, the algebra
$R(\beta)$ has a special quotient $R^{\Lambda}(\beta)$ 
which is called the {\it cyclotomic
Khovanov-Lauda-Rouquier algebra}. In \cite{KL09}, Khovanov and Lauda
conjectured that $\bigoplus_{\beta \in
Q^{+}}K\bl\Proj(\RL(\beta))\br$ has a $U_\A(\g)$-module
structure  and that there exists a $U_\A(\g)$-module isomorphism
$$V_\A(\Lambda) \simeq \bigoplus_{\beta \in Q^{+}} K\bl\Proj(R^{\Lambda}(\beta))\br,$$ 
where $V_\A(\Lambda)$ denotes the $U_\A(\g)$-module with highest weight
$\Lambda$. 
After partial results of
Brundan and Stroppel (\cite{BS08}),
Brundan and Kleshchev (\cite{BK08,BK09})
and Lauda and Vazirani (\cite{LV09}), 
the conjecture was proved by
Seok-Jin Kang and the author
for {\em all} symmetrizable Kac-Moody algebras (\cite{KK}).

For each $i \in I$, let us consider the restriction functor
and the induction functor:
\begin{equation*}
\begin{aligned}
& E_{i}^{\Lambda}\cl \Mod(R^{\Lambda}(\beta+\alpha_i)) \longrightarrow \Mod(R^{\Lambda}(\beta)), \\
& F_{i}^{\Lambda}\cl \Mod(R^{\Lambda}(\beta)) \longrightarrow
\Mod(R^{\Lambda}(\beta+ \alpha_i))
\end{aligned}
\end{equation*}
defined by
\begin{equation*}
\begin{aligned}
& E_{i}^{\Lambda}(N) = e(\beta, i) N = e(\beta, i)R^{\Lambda}(\beta+
\alpha_i) \otimes_{R^{\Lambda}(\beta+\alpha_i)} N, \\
& F_{i}^{\Lambda}(M) = R^{\Lambda}(\beta+\alpha_i) e(\beta, i)
\otimes_{R^{\Lambda}(\beta)} M,
\end{aligned}
\end{equation*}
where $M \in \Mod(R^{\Lambda}(\beta))$, $N \in
\Mod(R^{\Lambda}(\beta+\alpha_i))$.
Then these functors categorify the root operators $e_i$ and $f_i$
in the quantum groups.

It is obvious that $\E$ is a right adjoint functor of $\F$.

Khovanov-Lauda (\cite{KL09, KL08, KL10, L08}) and  Rouquier (\cite{R08})
conjectured that
$\E$ and $\F$ are biadjoint to each other.
Namely $\E$ is also a left adjoint of $\F$.
Furthermore they gave a candidate
of this adjunction explicitly from the first adjunction.
In this paper we prove that
this candidate gives indeed  adjunction for all cyclotomic \KLRs.

\bigskip
In order to prove this we use
a similar method employed in \cite{KK}.
Namely we use
the module $e(\beta,i^2)R(\beta+2\alpha_i)e(\beta+\alpha_i,i)
\otimes_{R(\beta+\al_i)}R^\Lambda(\beta+\alpha_i)$
in order to study $e(\beta,i^2)\RL(\beta+2\alpha_i)e(\beta+\alpha_i,i)$.
We fully use the fact that 
this module is a free right module
over the ring $\cor[x_{n+2}]$ (Lemma~\ref{lem:Grin}).

\smallskip
We mention that \cite{CL} and \cite{Web10} are related to our results.

\bigskip
This paper is organized as follows. In Section 2, we
recall the notions of \KLRs.
In Section 3, we recall the definition of
cyclotomic \KLRs\ and the results in \cite{KK}, and then
state our main result (Theorem~\ref{th:main}).
In Section 4, we interpret it in terms of the algebras
(\eqref{eq1}, \eqref{eq2} and \eqref{eq3}),
and we gave their proof in Section 5.
\vskip 5mm

{\it Acknowledgements.} 
We would like to thank 
Aaron Lauda by explaining his results with 
M.\ Khovanov, and also his recent paper \cite{CL} with S.\ Cautis.


\section{The Khovanov-Lauda-Rouquier algebra} \label{sec:R}
\subsection{Cartan data}
Let $I$ be a finite index set. An integral square matrix
$A=(a_{ij})_{i,j \in I}$ is called a {\em symmetrizable generalized
Cartan matrix} if it satisfies (i) $a_{ii} = 2$ $(i \in I)$, (ii)
$a_{ij} \le 0$ $(i \neq j)$, (iii) $a_{ij}=0$ if $a_{ji}=0$ $(i,j \in I)$,
(iv) there is a diagonal matrix
$D=\text{diag} (d_i \in \Z_{> 0} \mid i \in I)$ such that $DA$ is
symmetric.

A \emph{Cartan datum} $(A,P, \Pi,P^{\vee},\Pi^{\vee})$ consists of
\begin{enumerate}[(1)]
\item a symmetrizable generalized Cartan matrix $A$,
\item a free abelian group $P$ of finite rank, called the \emph{weight lattice},
\item $P^{\vee}\seteq\Hom(P, \Z)$, called the \emph{co-weight lattice},
\item $\Pi= \set{\alpha_i }{i \in I}\subset P$, called
the set of \emph{simple roots},
\item $\Pi^{\vee}= \set{h_i}{i \in I}\subset P^{\vee}$, called
the set of \emph{simple coroots},
\end{enumerate}
satisfying the condition:
$\langle h_i,\alpha_j \rangle = a_{ij}$ for all $i,j \in I$.

We denote by
$$P^{+} \seteq \set{ \lambda \in P}%
{\text{$\langle h_i, \lambda \rangle \in\Z_{\ge 0}$ for all $i \in I$}}$$
the set of \emph{dominant integral weights}.
The free abelian group $Q\seteq\soplus_{i \in I} \Z
\alpha_i$ is called the \emph{root lattice}. Set $Q^{+}= \sum_{i \in
I} \Z_{\ge 0} \alpha_i$. For $\alpha = \sum k_i \alpha_i \in Q^{+}$,
we define the {\it height}  $\haut(\alpha)$ of $\alpha$
to be $\haut(\alpha)=\sum k_i$.
Let $\mathfrak{h} = \Q \otimes_\Z P^{\vee}$. Since $A$ is
symmetrizable, there is a symmetric bilinear form $(\quad|\quad)$ on
$\mathfrak{h}^*$ satisfying
$$ (\alpha_i | \alpha_j) =d_i a_{ij} \quad (i,j \in I)
\quad\text{and $\lan h_i,\lambda\ran=
\dfrac{2(\alpha_i|\lambda)}{(\alpha_i|\alpha_i)}$ for any $\lambda\in\mathfrak{h}^*$ and $i \in I$}.$$

\subsection{Definition of  Khovanov-Lauda-Rouquier algebra} 

Let $(A, P, \Pi, P^{\vee}, \Pi^{\vee})$ be a Cartan datum. In this
section, we recall the construction of \KLR\ 
associated with $(A, P, \Pi, P^{\vee}, \Pi^{\vee})$
and its properties.
We take  as a base ring a graded commutative ring 
$\cor=\soplus\nolimits_{n\in\Z}\,\cor_n$ such that $\cor_n=0$ for any $n<0$.
Let us take a matrix $(Q_{ij})_{i,j\in I}$ in $\cor[u,v]$
such that $Q_{ij}(u,v)=Q_{ji}(v,u)$ and $Q_{ij}(u,v)$ has the form
\begin{equation} \label{eq:Q}
Q_{ij}(u,v) = \begin{cases}\hs{5ex} 0 \ \ & \text{if $i=j$,} \\
\sum\limits_{p,q\ge0}
t_{i,j;p,q} u^p v^q\quad& \text{if $i \neq j$,}
\end{cases}
\end{equation}
where $t_{i,j;p,q}\in \cor_{ -2(\alpha_i | \alpha_j)-(\alpha_i|\alpha_i) p - (\alpha_j|\alpha_j)q }$ and
 $t_{i,j}\seteq t_{i,j;-a_{ij},0} \in \cor_0^\times$.
In particular, we have $t_{i,j;p,q}=0$ if 
$(\alpha_i|\alpha_i) p + (\alpha_j|\alpha_j)q >-2(\alpha_i | \alpha_j)$.
Note that $t_{i,j;p,q} = t_{j,i;q,p}$.

We denote by
$S_{n} = \langle s_1, \ldots, s_{n-1} \rangle$ the symmetric group
on $n$ letters, where $s_i = (i, i+1)$ is the transposition.
Then $S_n$ acts on $I^n$.

\begin{Def}[\cite{KL09,{R08}}] \label{def:KLRalg}
The {\em \KLR}\ $R(n)$ of degree $n$
associated with a Cartan datum $(A, P, \Pi, P^{\vee}, \Pi^{\vee})$ and
$(Q_{ij})_{i,j\in I}$ is the associative algebra over $\cor$
generated by $e(\nu)$ $(\nu \in I^{n})$, $x_k$ $(1 \le k \le n)$,
$\tau_l$ $(1 \le l \le n-1)$ satisfying the following defining
relations:
\begin{equation} \label{eq:KLR}
\begin{aligned}
& e(\nu) e(\nu') = \delta_{\nu, \nu'} e(\nu), \ \
\sum_{\nu \in I^{n}}  e(\nu) = 1, \\
& x_{k} x_{l} = x_{l} x_{k}, \ \ x_{k} e(\nu) = e(\nu) x_{k}, \\
& \tau_{l} e(\nu) = e(s_{l}(\nu)) \tau_{l}, \ \ \tau_{k} \tau_{l} =
\tau_{l} \tau_{k} \ \ \text{if} \ |k-l|>1, \\
& \tau_{k}^2 e(\nu) = Q_{\nu_{k}, \nu_{k+1}} (x_{k}, x_{k+1})
e(\nu), \\
& (\tau_{k} x_{l} - x_{s_k(l)} \tau_{k}) e(\nu) = \begin{cases}
-e(\nu) \ \ & \text{if} \ l=k, \nu_{k} = \nu_{k+1}, \\
e(\nu) \ \ & \text{if} \ l=k+1, \nu_{k}=\nu_{k+1}, \\
0 \ \ & \text{otherwise},
\end{cases} \\[.5ex]
& (\tau_{k+1} \tau_{k} \tau_{k+1}-\tau_{k} \tau_{k+1} \tau_{k}) e(\nu)\\
&\hs{8ex} =\begin{cases} \dfrac{Q_{\nu_{k}, \nu_{k+1}}(x_{k},
x_{k+1}) - Q_{\nu_{k}, \nu_{k+1}}(x_{k+2}, x_{k+1})} {x_{k} -
x_{k+2}}e(\nu) \ \ & \text{if} \
\nu_{k} = \nu_{k+2}, \\
0 \ \ & \text{otherwise}.
\end{cases}
\end{aligned}
\end{equation}

%
\end{Def}
\noindent
Note that $R(n)$ has an anti-involution $\psi$ that fixes the
generators $x_k$, $\tau_l$ and $e(\nu)$.

The $\Z$-grading on $R(n)$ is given by
\begin{equation} \label{eq:Z-grading}
\deg e(\nu) =0, \quad \deg\; x_{k} e(\nu) = (\alpha_{\nu_k}
| \alpha_{\nu_k}), \quad\deg\; \tau_{l} e(\nu) = -
(\alpha_{\nu_l} | \alpha_{\nu_{l+1}}).
\end{equation}

For $a,b,c \in \{1,\ldots,n\}$, we define the elements of $R(n)$ by
\begin{equation}
\begin{aligned}
& e_{a,b} = \sum_{\nu \in I^{n},\, \nu_{a}=\nu_{b}} e(\nu), \\
& Q_{a,b} = \sum_{\nu \in I^n} Q_{\nu_{a},\, \nu_{b}} (x_{a},
x_{b}) e(\nu), \\
& \overline{Q}_{a,b,c} = \sum_{\nu \in I^{n},\; \nu_{a} = \nu_{c}}
\dfrac{Q_{\nu_{a}, \nu_{b}}(x_{a}, x_{b}) - Q_{\nu_{a},
\nu_{b}}(x_{c}, x_{b})}{x_{a} - x_{c}} e(\nu) \quad\text{if $a\not=c$.}
\end{aligned}\label{def:Q}
\end{equation}
Then we have
\begin{equation}
\begin{aligned}
& Q_{a,b} = Q_{b,a}, \quad \tau_{a}^2 = Q_{a, a+1}, \\
& \tau_{a+1} \tau_{a} \tau_{a+1} = \tau_{a} \tau_{a+1} \tau_{a} +
\overline{Q}_{a, a+1, a+2}.
\end{aligned}\label{eq:tau3}
\end{equation}
We define the operators $\partial_{a,b}$ on
$\soplus_{\nu\in I^n}\cor[x_1, \ldots, x_n]e(\nu)$.
by
\begin{equation}
\partial_{a,b} f = \dfrac{s_{a,b} f - f} {x_{a} - x_{b}}e_{a,b}, \quad
\partial_{a} = \partial_{a,a+1},
\end{equation}
where $s_{a,b} = (a,b)$ is the transposition.

Thus we obtain
\begin{equation} \label{eq:partial}
\begin{aligned}
& \tau_{a} e_{b,c} = e_{s_a(b), s_a(c)} \tau_{a}, \\
& \tau_{a} f - (s_a f) \tau_{a} = f \tau_{a} - \tau_{a} (s_a f) =
(\partial_{a} f) e_{a,a+1}.
\end{aligned}
\end{equation}

For $n\in \Z_{\ge 0}$ and $\beta \in Q^{+}$ such that $\haut(\beta)=n$, we set
$$I^{\beta} = \set{ \nu = (\nu_1, \ldots, \nu_n) \in I^n }%
{\alpha_{\nu_1} + \cdots + \alpha_{\nu_n} = \beta }.$$
We define
\eq
&&\ba{l}
e(\beta) = \sum_{\nu \in I^{\beta}} e(\nu), \\[1ex]
 R(\beta) = R(n) e(\beta)=\soplus_{\nu\in I^\beta}R(n)e(\nu).
\ea\eneq

The algebra $R(\beta)$ is called the {\it {\KLR} at $\beta$}.

For $\ell\ge0$, we set
\eq
&&\ba{l}
e(\beta, i^\ell) = \sum_{\nu}e(\nu)\in R(\beta+\ell\alpha_i)\\[1ex]
\hs{10ex}\parbox{55ex}{where $\nu$ ranges over the set of
$\nu\in I^{\beta+\ell\al_i}$ such that
$\nu_k=i$ for $n+1\le k\le n+\ell$.} 
\ea\eneq

We sometimes regard $R(\beta)$ as a $\cor$-subalgebra of 
the $\cor$-algebra $e(\beta,i^\ell)R(\beta+\ell\al_i)e(\beta,i^\ell)$.


%

\Th\label{th:R(n)xR(n)}
Let $\beta\in Q^+$ with $\haut(\beta)=n$ and $i\in I$.
Then there exists a natural isomorphism
\eq
&&
\ba{l}
R(\beta) e(\beta-\al_i,i) \otimes_{R(\beta-\al_i)}\cor\tau_n\tens e(\beta-\al_i, i) R(\beta)
\soplus\cor[x_{n+1}]\otimes R(\beta)\\[1ex]
\hs{5ex}\isoto e(\beta,i) R(\beta+\al_i) e(\beta,i).
\ea
\eneq
Here $R(\beta) e(\beta-\al_i,i) \otimes_{R(\beta-\al_i)}\cor\tau_n\tens e(\beta-\al_i, i) R(\beta)\to e(\beta,i) R(\beta+\al_i) e(\beta,i)$
is given by $a\otimes \tau_{n}\tens b\mapsto a\tau_{n} b$.
\entheorem
Here, $\tau_n$ in $\cor \tau_n$ is
a symbolical basis of a free $\cor$-module of rank one.
We sometimes use such notations in order to
make morphisms more explicit.

Note that if $\beta-\al_i\not\in Q^+$ then
$R(\beta) e(\beta-\al_i,i) \otimes_{R(\beta-\al_i)}
\cor\tau_n\tens e(\beta-\al_i, i) R(\beta)$ should be understood
to be zero.
\section{The cyclotomic Khovanov-Lauda-Rouquier algebras} \label{sec:RLambda}

\subsection{Definition of cyclotomic \KLRs}
Let $\Lambda\in P^+$ be a dominant integral weight.
For each $i\in I$, we shall choose a monic polynomial of degree 
$\langle h_{i},\Lambda \rangle$
\eq
a_i^\Lambda(u)=\sum_{k=0}^{\langle h_{i},\Lambda \rangle}
c_{i;k}u^{\langle h_{i},\Lambda \rangle-k}
\eneq
with $c_{i;k}\in \cor_{k(\alpha_i|\alpha_i)}$ and $c_{i;0}=1$.

For $k$ ($1\le k\le n$) and $\beta\in Q^+$ with $\haut(\beta)=n$, we set
\begin{equation}
\x[k]= \sum_{\nu \in I^{\beta}} a_{\nu_k}^\Lambda(x_k) e(\nu)\in R(\beta).
\end{equation}
Hence $\x[k]e(\nu)$ is a homogeneous element of $R(\beta)$ with degree 
$2(\alpha_{\nu_k}|\Lambda)$.

\begin{Def} \label{def:RLambda} For $\beta\in Q^+$
the {\em cyclotomic Khovanov-Lauda-Rouquier
algebra $R^{\Lambda}(\beta)$ at $\beta$} is defined to be the
quotient algebra
$$R^{\Lambda}(\beta) = \dfrac{R(\beta)}  {R(\beta)\x R(\beta)}.$$
\end{Def}

In this paper we {\em forget the grading},
and we denote by $\Mod(\RL(\beta))$ the abelian
category of $\RL(\beta)$-modules.

For each $i\in I$, we define the functors
\begin{equation*}
\begin{aligned}
& \E\cl \Mod(R^{\Lambda}(\beta+\alpha_i)) \longrightarrow \Mod(R^{\Lambda}(\beta)), \\
& \F\cl \Mod(R^{\Lambda}(\beta)) \longrightarrow
\Mod(R^{\Lambda}(\beta+ \alpha_i))
\end{aligned}
\end{equation*}
by
\begin{equation} \label{eq:E_iLambda}
\begin{aligned}
&\E(N) = e(\beta, i) N \simeq e(\beta, i)R^{\Lambda}(\beta+
\alpha_i) \otimes_{R^{\Lambda}(\beta+\alpha_i)} N\\
&\hs{17.5ex}\simeq\Hom_{\RL(\beta+\al_i)}(\RL(\beta+\al_i)e(\beta,i),N), \\
&\F(M) = R^{\Lambda}(\beta+\alpha_i) e(\beta, i)
\otimes_{R^{\Lambda}(\beta)} M,
\end{aligned}
\end{equation}
where $M \in \Mod(R^{\Lambda}(\beta))$ and $N \in
\Mod(R^{\Lambda}(\beta+\alpha_i))$.

Then the following result is proved in \cite{KK}.

\begin{Thm}[\cite{KK}]\label{th:proj}
The module $R^{\Lambda}(\beta+\alpha_i)e(\beta, i)$
is a projective right $R^\Lambda(\beta)$-module.
Similarly, $e(\beta, i)R^{\Lambda}(\beta+\alpha_i)$
is a projective left $R^\Lambda(\beta)$-module.
 \end{Thm}

\Cor\label{cor:exact} \hfill
\bnum
\item
The functor $\E$ sends finitely generated
projective modules to finitely generated projective modules.
\item The functor $\F$ is exact.
\enum
\encor

%
%

\subsection{}
The pair $(\F,\E)$ has a canonical adjunction
: the unit $\eta\cl\id\To \E\F$ 
and the counit $\eps\cl\F\E\To\id$.

For $\beta\in Q^+$ with $\haut(\beta)=n$, the functors 
$$\xymatrix@C=7ex{\Mod(R^\La(\beta))\ar@<.5ex>[r]^-{\F}
&\Mod(R^\La(\beta+\al_i))\ar@<.5ex>[l]^-{\E}}$$
are represented by the kernel bimodules
$R^\La(\beta+\al_i)e(\beta,i)$ and $e(\beta,i)R^\la(\beta+\al_i)$
as in \eqref{eq:E_iLambda}.
In the sequel, we denote by $\One_\beta$ the identity functor of
the category $\Mod(\RL(\beta))$, and
we denote by $\One_\beta\E=\E\One_{\beta+\al_i}$
the restriction functor $\E\cl\Mod(\RL(\beta+\al_i))\to\Mod(\RL(\beta))$.
Similarly, $\F\One_\beta=\One_{\beta+\al_i}\F$ denotes
the induction functor $\F\cl\Mod(\RL(\beta))\to\Mod(\RL(\beta+\al_i))$.

Let us denote by $x$ the endomorphism of $\One_{\beta}\E$ 
represented by the left multiplication
of $x_{n+1}$ on $e(\beta,i)R^\La(\beta+\al_i)$
and by $\tau$ the endomorphism of 
$\One_{\beta}\E\E\cl\Mod(\RL(\beta+2\al_i))\to\Mod(\RL(\beta))$
represented by the left multiplication of
$\tau_{n+1}$ on $e(\beta,i)R^\La(\beta+\al_i)
\otimes_{R^\La(\beta+\al_i)}e(\beta+\al_i,i)R^\La(\beta+2\al_i)
\simeq e(\beta,i^2)R^\La(\beta+2\al_i)$.
Similarly the endomorphism $x$ of $\F\One_{\beta}$
is represented by the right multiplications
of $x_{n+1}$ on $R^\La(\beta+\al_i)e(\beta,i)$
and the endomorphism $\tau$ of $\F\F\One_{\beta}
\cl\Mod(\RL(\beta))\to\Mod(\RL(\beta+2\al_i))$
is represented by the right multiplication of
$\tau_{n+1}$ on $R^\La(\beta+2\al_i)e(\beta+\al_i,i)
\otimes_{R^\La(\beta+\al_i)}R^\La(\beta+\al_i)e(\beta,i)
\simeq R^\La(\beta+2\al_i)e(\beta,i^2)$.
Then $x\in\End(\F\One_{\beta})$ and $x\in\End(\One_{\beta}\E)$ 
are dual to each other
and  $\tau\in\End(\F\F\One_{\beta})$ and $\tau\in\End(\One_{\beta}\E\E)$ are 
dual to each other.

By the adjunction, $\tau\in \End(\E\E)$ induces a morphism
\eq
\sigma\cl\F\E\One_\beta\To\E\F\One_\beta.
\eneq
It is represented by the morphism
$$\RL(\beta)e(\beta-\al_i,i)\otimes_{\RL(\beta-\al_i)}e(\beta-\al_i,i)\RL(\beta)
\to e(\beta,i)\RL(\beta+\al_i)e(\beta,i)$$
given by $x\otimes y\To x\tau_ny$.

The following theorem was formulated as one of
the axioms for the categorification of representations of
quantum groups (\cite{CR08,KL10,L08,R08}),
and proved in \cite{KK} for an arbitrary \KLR.
\begin{Thm}[\cite{KK}] \label{thm:M}
Set $\lambda \seteq \Lambda - \beta$ and $\la_i\seteq\lan h_i,\la\ran$.
\bna
\item Assume $\la_i\seteq\langle h_{i}, \lambda \rangle \ge 0$.
The the morphism of endofunctors on $\Mod(R^{\Lambda}(\beta))$
\eqn
\rho&:&\F\E\One_\beta
\oplus \soplus_{k=0}^{\la_i-1}\cor x^k\otimes \One_\beta
\To \E\F\One_\beta
\eneqn 
is an isomorphism.
Here $\F\E\One_\beta\to\E\F\One_\beta$ is given by $\sigma$,
and $\cor x^k\otimes \One_\beta\to\F\E\One_\beta$ is given by
$(x^k\F)\circ\eta=(\E x^k)\circ\eta\cl\One_\beta\to \E\F\One_\beta$.
\item Assume that $\la_i\le 0$.
Then the morphism 
\eqn
\rho &:&\F\E\One_\beta\To \E\F\One_\beta
\oplus \soplus_{k=0}^{-\la_i-1}\cor (x^{-1})^k\otimes \One_\beta
\eneqn
is an isomorphisms. Here $\F\E\One_\beta\to\E\F\One_\beta$ is given by $\sigma$,
and $\F\E\One_\beta\to \cor (x^{-1})^k\otimes \One_\beta$
is given by
$\eps\circ(x^k\E)=\eps\circ(\F x^k)\cl\F\E\One_\beta\to \One_\beta$.
\end{enumerate}
\end{Thm}
In the theorem, $x^k$ in $\cor x^k$ and $(x^{-1})^k$ in $\cor (x^{-1})^k$ are
a symbolical basis of a free $\cor$-module.

Now let us define the morphism
$\heta\cl\One_\beta\to\F\E\One_\beta$ as follows.
\bnum
\item
If $\la_i
\seteq\lan h_i,\la\ran\ge0$, then $\heta$ is given by the commutativity of
\eqn
&&\xymatrix@C=9ex{
{\F\E\One_\beta}&{\F\E\One_\beta
\oplus \soplus_{k=0}^{\la_i-1}\cor x^k\otimes \One_\beta}
\ar[l]_-{\text{projection}}\ar[d]_-{\bigwr}^-{\rho}
\\
\One_\beta\ar[r]^{x^{\la_i}F\circ\eta}
\ar[u]_{-\heta}
&\E\F.
}
\eneqn
Here the top horizontal arrow is the projection.
The minus sign in front of $\heta$ should be noted.
\item
If $\la_i<0$, then $\widehat{\eta}$ is defined as
the composition
\eqn
&&\xymatrix@C=7ex{
{\One_\beta}\ar[r]^-{\widehat{\eta}}\ar[d]_-{\bigwr}
&{\F\E\One_\beta}\ar[d]_-{\bigwr}^-{\rho}\\
\cor (x^{-1})^{-\la_i-1}\otimes \One_\beta\ar@{^{(}->}[r]&
\E\F\One_\beta
\oplus \soplus_{k=0}^{-\la_i-1}\cor (x^{-1})^k\otimes \One_\beta.
}
\eneqn
Here the bottom horizontal arrow is the canonical inclusion and
the left vertical arrow is derived from
$\cor\isoto \cor (x^{-1})^{-\la_i-1}$ (
$1\mapsto (x^{-1})^{-\la_i-1}$).
\ee

The morphism
$\widehat{\eps}\cl\E\F\One_\beta\to\One_\beta$ is defined as follows.
\bnum
\item
If $\la_i>0$, then $\widehat{\eps}$ is defined as
the composition
\eqn
&&\xymatrix@C=9ex{
{\F\E\One_\beta
\oplus \soplus_{k=0}^{\la_i-1}\cor x^k\otimes \One_\beta}
\ar[r]^-{\text{projection}}\ar[d]_-{\bigwr}^-{\rho}&
{\cor x^{\la_i-1}\otimes \One_\beta}
\ar[d]_-{\bigwr}\\
\E\F\One_\beta\ar[r]^{\widehat{\eps}}&
\One_\beta.
}
\eneqn
Here the top horizontal arrow is the canonical projection and
the right vertical arrow is induced by $x^{\la_i-1}\mapsto 1$.

\item
If $\la_i\le0$, then $\heps$ is defined as
the composition
\eqn
&&\xymatrix@C=7ex{
\One_\beta&{\F\E\One_\beta}\ar[d]_-{\bigwr}^-{\rho}%
\ar[l]_-{\eps\circ(x^{-\la_i}\E)}
\\
\E\F
\ar[u]_{\heps}\ar@{^{(}->}[r]
&
\E\F\One_\beta
\oplus \soplus_{k=0}^{-\la_i-1}\cor (x^{-1})^k\otimes \One_\beta.
}
\eneqn
Here the bottom horizontal arrow is the canonical inclusion.
\ee

Now our main result can be stated as follows.
\Th\label{th:main}
The pair $(\E,\F)$ is an adjoint pair
with $(\heta,\heps)$ as adjunction.
Namely the compositions $\E\To[\E\,\heta] \E\F\E\To[\heps\;\E] \E$
and $\F\To[\heta\,\F] \F\E\F\To[\F\,\heps]\F$ are equal to the identities.
\entheorem
We shall prove this theorem in the rest of the paper.

\section{Proof of Theorem \ref{th:main}}
\subsection{}
We shall first prove that the composition
$\One_\beta\E\To[\E\,\heta] \One_\beta\E\F\E\To[\heps\;\E] \One_\beta\E$
is equal to the identity.
Here $\beta\in Q^+$ with $\haut(\beta)=n$ and 
we set $\la\seteq\La-\beta$ and $\la_i\seteq\lan h_i,\la\ran$.

\subsubsection{$\la_i\ge2$ Case}\quad
We shall first assume that
$\la_i\ge2$.
Then the composition
$\One_\beta\E\To[\E\,\heta] \One_\beta\E\F\E\To[\heps\;\E] \One_\beta\E$
can be described by the kernel bimodules as follows.
The morphism $\One_\beta\E\To[\E\,\heta] \One_\beta\E(\F\E\One_{\beta+\al_i})$ 
is given by the $\bl\RL(\beta),\RL(\beta+\al_i)\br$-bilinear homomorphism:
\eqn
&&\xymatrix{
e(\beta,i)\RL(\beta+\al_i)\ar[d]_{-x_{n+2}^{\la_i-2}}\\
e(\beta,i^2)\RL(\beta+2\al_i)e(\beta+\al_i,i)\\
\txt{\parbox{55ex}
{$e(\beta,i)\RL(\beta+\al_i)e(\beta,i)\tens\limits_{\RL(\beta)}
\cor\tau_{n+1}\tens e(\beta,i)\RL(\beta+\al_i)$\\
\hs{4ex}\hfill$\oplus\kern-1ex\soplus_{k=0}^{k=\la_i-3}
\cor x_{n+2}^k\otimes e(\beta,i)\RL(\beta+\al_i)
$}
}\ar[u]^-{\rho}_-\bigwr\ar[d]^-{\txt{projection}}\\
e(\beta,i)\RL(\beta+\al_i)e(\beta,i)\tens\limits_{\RL(\beta)}
\cor\tau_{n+1}\tens e(\beta,i)\RL(\beta+\al_i).
}
\eneqn

The morphism $(\One_\beta\E\F)\E\To[\heps\;\E] \One_\beta\E$
is given by the $\bl\RL(\beta),\RL(\beta+\al_i)\br$-bilinear homomorphism:
\eqn
&&\xymatrix{
e(\beta,i)\RL(\beta+\al_i)e(\beta,i)\tens\limits_{\RL(\beta)}
\cor\tau_{n+1}\tens e(\beta,i)\RL(\beta+\al_i)\\
\txt{\parbox{61ex}{
$\Bigl(\RL(\beta)e(\beta-\al_i,i)\kern-2ex\tens\limits_{\RL(\beta-\al_i)}\kern-2ex
\cor\tau_n\tens e(\beta-\al_i,i)\RL(\beta)$\\[1ex]
\hs{5ex}$\oplus\kern-0ex\soplus_{k=0}^{\la_i-1}
\cor x_{n+1}^k\otimes\RL(\beta)\Bigr)
\tens_{\RL(\beta)}\cor\tau_{n+1}\tens e(\beta,i)\RL(\beta+\al_i)
$}}\ar[u]^-{\rho}_-{\bigwr}\ar[d]^-{\text{projection}}\\
\cor x_{n+1}^{\la_i-1}\tens \cor\tau_{n+1}\tens 
e(\beta,i)\RL(\beta+\al_i)
\ar[d]^-{\bigwr}\\
e(\beta,i)\RL(\beta+\al_i).
}
\eneqn
Hence in order to see that the composition is the identity,
it is enough to show the inclusion
\eq
&&\ba{l}
x_{n+2}^{\la_i-2}e(\beta,i^2)+x_{n+1}^{\la_i-1}\tau_{n+1}e(\beta,i^2)\\[1ex]
\hs{5ex}\in 
\RL(\beta)\tau_n\tau_{n+1}e(\beta-\al_i,i^3)\RL(\beta+\al_i)\\[1ex]
\hs{10ex}+\sum_{k=0}^{\la_i-2}x_{n+1}^k\tau_{n+1}e(\beta,i^2)
\RL(\beta+\al_i)\\[1ex]
\hs{15ex}+\sum_{k=0}^{\la_i-3}x_{n+2}^ke(\beta,i^2)\RL(\beta+\al_i)
\ea\label{eq1}
\eneq
as an element of $e(\beta,i^2)\RL(\beta+2\al_i)e(\beta+\al_i,i)$.

This inclusion is proved in \S\,\ref{sec:proof}.

\subsubsection{}
Now let us treat the case $\la_i=1$.

The morphism $\One_\beta\E\To[\E\,\heta] \One_\beta\E(\F\E\One_{\beta+\al_i})$ 
is given by 
\eqn
&&\xymatrix{
e(\beta,i)\RL(\beta+\al_i)\ar[d]^{\text{inclusion}}\\
e(\beta,i^2)\RL(\beta+2\al_i)e(\beta+\al_i,i)
\oplus e(\beta,i)\RL(\beta+\al_i) \\
e(\beta,i)\RL(\beta+\al_i)e(\beta,i)\tens_{\RL(\beta)}
e(\beta,i)\RL(\beta+\al_i)\ni u.\ar[u]^{\bigwr}_{\rho=\varSigma\oplus\tE}
}
\eneqn
(See below for $\varSigma$ and $\tE$.)

The morphism $(\One_\beta\E\F)\E\To[\heps\;\E] \One_\beta\E$
is given by

\eqn
&&\xymatrix{
e(\beta,i)\RL(\beta+\al_i)e(\beta,i)\tens_{\RL(\beta)}
e(\beta,i)\RL(\beta+\al_i)\ni u\\
\txt{\parbox{55ex}{
$\Bigl(\RL(\beta)e(\beta-\al_i,i)\kern-2ex\tens\limits_{\RL(\beta-\al_i)}\kern-2ex
e(\beta-\al_i,i)\RL(\beta)\oplus\RL(\beta)\Bigr)$\\[1ex]
\hs{1ex}\hfill$
\tens_{\RL(\beta)}
e(\beta,i)\RL(\beta+\al_i)$}}
\ar[u]_-{\rho}^-{\bigwr}\ar[d]^-{\text{projection}}\\
e(\beta,i)\RL(\beta+\al_i).
}
\eneqn
Hence in order to see that the composition is the identity,
it is enough to show the following existence :
\eq\label{eq2}
&&
\parbox{73ex}{
There exists $u\in 
e(\beta,i)\RL(\beta+\al_i)e(\beta,i)\tens_{\RL(\beta)}
e(\beta,i)\RL(\beta+\al_i)$ such that

\hs{3ex}
\parbox{67ex}{
\bna
\item
 $\varSigma(u)=0$,
\item $\tE(u)=e(\beta,i)$,
\item $u-e(\beta,i)\tens e(\beta,i)
\in \bl \RL(\beta) e(\beta-\al_i,i^2)\tau_ne(\beta-\al_i,i^2)\RL(\beta)\br$\hfill\\
\hs{30ex}\hfill$\tens_{\RL(\beta)}e(\beta,i)\RL(\beta+\al_i)$.
\ee}
}
\eneq
Here 
\eqn
\varSigma&\cl& e(\beta,i)\RL(\beta+\al_i)e(\beta,i)\tens_{\RL(\beta)}
e(\beta,i)\RL(\beta+\al_i)\\
&&\hs{30ex}\To e(\beta,i^2)\RL(\beta+2\al_i)e(\beta+\al_i,i)\nonumber
\eneqn
is given by $\Sigma(a\otimes b)=a\tau_{n+1}b$,
and 
\eqn
\tE&\cl& e(\beta,i)\RL(\beta+\al_i)e(\beta,i)\tens_{\RL(\beta)}
e(\beta,i)\RL(\beta+\al_i)\to e(\beta,i)\RL(\beta+\al_i)
\eneqn
is given by $\tE(a\otimes b)=ab$.

The proof of \eqref{eq2} will be given in \S\,\ref{sec:proof}.

\subsubsection{}

Now we assume that $\la_i\le0$.
Then the composition
$\One_\beta\E\To[\E\,\heta] \One_\beta\E\F\E\To[\heps\;\E] \One_\beta\E$
can be described by the kernel bimodules as follows.

The morphism 
$\One_\beta\E\To[\E\,\heta] \One_\beta\E(\F\E\One_{\beta+\al_i})$ is given by 
\eqn
&&\xymatrix{
e(\beta,i)\RL(\beta+\al_i)\ar[d]^{\bigwr}\\
\cor (x_{n+1}^{-1})^{1-\la_i}\tens e(\beta,i)\RL(\beta+\al_i)
\ar[d]^-{\txt{inclusion}}\\
\txt{\parbox{61ex}
{$
e(\beta,i^2)\RL(\beta+2\al_i)e(\beta+\al_i,i)
\oplus\kern-.5ex\soplus_{k=0}^{1-\la_i}
\cor (x_{n+1}^{-1})^k\otimes e(\beta,i)\RL(\beta+\al_i)
$}}
\\
e(\beta,i)\RL(\beta+\al_i)e(\beta,i)\tens_{\RL(\beta)}
e(\beta,i)\RL(\beta+\al_i).
\ar[u]^-{\rho=f\oplus \oplus_kH'_k}_-{\bigwr}
}
\eneqn

The morphism $(\One_\beta\E\F)\E\To[\heps\;\E] \One_\beta\E$
is given by
\eqn
&&\xymatrix{
e(\beta,i)\RL(\beta+\al_i)e(\beta,i)\tens_{\RL(\beta)}
e(\beta,i)\RL(\beta+\al_i)\ar[d]^-{\text{inclusion}}\\
\txt{\parbox{61ex}{
$e(\beta,i)\RL(\beta+\al_i)e(\beta,i)\tens\limits_{\RL(\beta)}
e(\beta,i)\RL(\beta+\al_i)$\\
\hs{10ex}$\oplus\kern-1ex\soplus_{k=0}^{-\la_i-1}
\cor (x_{n}^{-1})^k\otimes e(\beta,i)\RL(\beta+\al_i)$
}}\\
\RL(\beta)e(\beta-\al_i,i)\tens_{\RL(\beta-\al_i)}
e(\beta-\al_i,i^2)\RL(\beta+\al_i)\ni v
\ar[u]^-{\rho=g\oplus\mathop\oplus\limits_kT_k}_-{\bigwr}\ar[d]^-{T_{-\la_i}}\\
e(\beta,i)\RL(\beta+\al_i).
}
\eneqn
Hence in order to see that the composition is the identity,
it is enough to show the following:
\eq
&&\left\{\hs{.5ex}
\parbox{70ex}{
There exists
$v\in\RL(\beta)e(\beta-\al_i)\tens_{\RL(\beta-\al_i)}e(\beta,i)\RL(\beta+\al_i)$ such that we have

\hs{5ex}\parbox{40ex}{
\bna
\item $T_k(v)=0$ for $0\le k\le-\la_i-1$
\item $T_{-\la_i}(v)=e(\beta,i)$,
\item $G(v)=0$,
\item $H_k(v)=0$ for $0\le k\le-\la_i$,
\item $H_{1-\la_i}(v)=e(\beta,i)$.
\ee
}
}\right.\label{eq3}
\eneq
Here the homomorphism
$$f\cl e(\beta,i)\RL(\beta+\al_i)e(\beta,i)\tens_{\RL(\beta)}
e(\beta,i)\RL(\beta+\al_i)\to e(\beta,i^2)\RL(\beta+2\al_i)e(\beta+\al_i,i)$$
is given by $f(a\otimes b)=a\tau_{n+1}b$, 
\eqn&&H'_k\cl e(\beta,i)\RL(\beta+\al_i)e(\beta,i)\tens_{\RL(\beta)}
e(\beta,i)\RL(\beta+\al_i)\\
&&\hs{10ex}\To\cor (x_{n+1}^{-1})^k\otimes e(\beta,i)\RL(\beta+\al_i)\simeq
 e(\beta,i)\RL(\beta+\al_i)
\eneqn
is given by 
$H'_k(a\otimes b)=ax_{n+1}^kb$, 
$$g\cl \RL(\beta)e(\beta-\al_i)\kern-2ex\tens_{\RL(\beta-\al_i)}\kern-2ex
e(\beta,i)\RL(\beta+\al_i)
\to e(\beta,i)\RL(\beta+\al_i)e(\beta,i)
\kern-1ex\tens\limits_{\RL(\beta)}\kern-1ex
e(\beta,i)\RL(\beta+\al_i)$$
is given by
$g(a\tens b)=a\tau_n\tens b$,
$$T_k\cl \RL(\beta)e(\beta-\al_i)
\tens_{\RL(\beta-\al_i)}e(\beta,i)\RL(\beta+\al_i)
\to e(\beta,i)\RL(\beta+\al_i)$$
is given by $T_k(a\tens b)=ax_{n}^kb$,
$$G=f\circ g\cl \RL(\beta)e(\beta-\al_i)
\tens_{\RL(\beta-\al_i)}e(\beta,i)\RL(\beta+\al_i)
\to e(\beta,i^2)\RL(\beta+2\al_i)e(\beta+\al_i,i)$$ 
is given by $G(a\otimes b)=a\tau_n\tau_{n+1}b$, and
$$H_k=H'_k\circ g\cl \RL(\beta)e(\beta-\al_i)
\tens_{\RL(\beta-\al_i)}e(\beta,i)\RL(\beta+\al_i)
\to e(\beta,i)\RL(\beta+\al_i)$$
is given by $H_k(a\otimes b)=a\tau_nx_{n+1}^kb$.

The statement \eqref{eq3} is proved in \S\,\ref{sec:proof}.

\subsection{}
Let us show that the composition 
$\F\To[\heta\,\F] \F\E\F\To[\F\,\heps]\F$
is equal to the identity by reducing it to the corresponding
statement for $\E\To[\E\,\heta] \E\F\E\To[\heps\;\E] \E$.

\bigskip
Let us recall that
$\psi$ is the anti-involution of $\RL(\beta)$
sending the generators $e(\nu)$, $x_k$, $\tau_k$ to themselves.
For an $\RL(\beta)$-module $M$,
we denote by $M^\psi$ the $\RL(\beta)^\op$-module induced by $\psi$ from $M$,
where $\RL(\beta)^\op$ is the opposite ring of $\RL(\beta)$.
We define the bifunctor
$$\Psi_\beta\cl\Mod(\RL(\beta))\times \Mod(\RL(\beta))\To\Mod(\cor)$$
by
$$\Psi_\beta(M,N)\seteq M^\psi\tens_{\RL(\beta)}N.$$
We have a functorial isomorphism 
$$\Psi_\beta(M,N)\simeq\Psi_\beta(N,M)
\quad\text{in $M,N\in \Mod(\RL(\beta))$.}$$

For two $\cor$-linear categories $\shc$ and $\shc'$,
let us denote by $\Fct_\cor(\shc,\shc')$ be the category of
$\cor$-linear functors from $\shc$ to $\shc'$.
Then $\Psi_\beta$ induces  a functor
$$\HH_\beta\cl\Mod(\RL(\beta))\To \Fct_\cor\bl\Mod(\RL(\beta)),\Mod(\cor)\br$$
by assigning to $M\in\Mod(\RL(\beta))$ the functor $N\mapsto\Psi_\beta(M,N)$.
The following lemma similar to Yoneda lemma
is easily proved, and its proof is omitted.
\Lemma\label{lem:Yoneda} The functor $\HH_\beta$ is fully faithful.
\enlemma
For $\beta,\beta'\in Q^+$
and a pair of $\cor$-linear functors $F\cl\Mod(\RL(\beta))\to\Mod(\RL(\beta'))$
and $G\cl\Mod(\RL(\beta'))\to\Mod(\RL(\beta))$,
we say that $F$ and $G$
are {\em $\Psi$-adjoint} or $G$ is a {\em $\Psi$-adjoint} of $F$ if
there exists a functorial isomorphism 
$$\Psi_{\beta'}(F(M),N)\simeq \Psi_\beta(M,G(N))\quad
\text{in $M\in\Mod(\RL(\beta))$ and $N\in\Mod(\RL(\beta'))$.}$$
For a given $F$, a $\Psi$-adjoint of $F$ is unique
up to a unique isomorphism if it exists.
We shall denote by $F^\vee$ the $\Psi$-adjoint of $F$
(if it exists).

If $\Mod(\RL(\beta))\To[F]\Mod(\RL(\beta'))\To[F']
\Mod(\RL(\beta''))$ are functors which admit $\Psi$-adjoint,
then $F^\vee\circ F'{}^\vee$ is a $\Psi$-adjoint of $F'\circ F$.

Now let $F_k\cl \Mod(\RL(\beta))\to\Mod(\RL(\beta'))$ ($k=1,2$)
be two functors.
Then Lemma~\ref{lem:Yoneda} implies
$$\Hom(F_1,F_2)\simeq \Hom(F_1^\vee,F_2^\vee).$$
For $f\in \Hom(F_1,F_2)$, the corresponding morphism
in $\Hom(F_1^\vee,F_2^\vee)$ is called the $\Psi$-adjoint of $f$
and we denote it by $f^\vee$. By the definition we have a commutative diagram
$$\xymatrix{
\Psi_{\beta'}(F_1(M),N)\ar@{-}[r]^\sim\ar[d]^f&\Psi_\beta(M,F_1^\vee(N))\ar[d]^{f^\vee}\\
\Psi_{\beta'}(F_2(M),N)\ar@{-}[r]^\sim&\Psi_\beta(M,F_2^\vee(N)).
}$$
Then $(f\circ g)^\vee=f^\vee\circ g^\vee$ for $F_1\To[g] F_2\To[f]F_3$.

The following lemma is elementary and its proof is omitted.
\Lemma
\bnum
\item Let $K$ be a $\bl\RL(\beta'),\RL(\beta)\br$-bimodule 
and the functor $F\cl \Mod(\RL(\beta))\to\Mod(\RL(\beta'))$
is given by $K\tens_{\RL(\beta)}\scbul$.
Then 
$F$ admits  a $\Psi$-adjoint.
\item Conversely if a $\cor$-linear functor
 $F\cl \Mod(\RL(\beta))\to\Mod(\RL(\beta'))$ admits a $\Psi$-adjoint, then
$F$ is isomorphic to
$F(\RL(\beta))\tens_{\RL(\beta)}\scbul$,
and $F^\vee(\RL(\beta'))\simeq F(\RL(\beta))^\psi$ as 
$\bl\RL(\beta),\RL(\beta')\br$-bimodules.
\ee
\enlemma

\bigskip
We can easily see that $\E$ and $\F$ are $\Psi$-adjoint. Moreover,
$x\in\End(\E)$ and $x\in\End(\F)$, $\tau\in\End(\E\circ\E)$
and $\tau\in\End(\F\circ\F)$ are $\Psi$-adjoint, respectively.
We can also see that $\eta\in\Hom(\One_\beta,\E\F\One_\beta)$
is a $\Psi$-adjoint of itself.
Similarly $\eps\in\Hom(\F\E\One_\beta,\One_\beta)$
$\sigma\in\Hom(\F\E,\E\F)$ are $\Psi$-adjoint of themselves.
Note that $\F\E$ and $\E\F$ are a $\Psi$-adjoint of themselves.
Hence $\heta$ and $\heps$ are also a $\Psi$-adjoint of themselves.

\smallskip
Therefore $\F\To[\heta\,\F] \F\E\F\To[\F\,\heps]\F$
is a $\Psi$-adjoint of $\E\To[\E\,\heta] \E\F\E\To[\heps\;\E] \E$.
Hence if the composition of $\E\To[\E\,\heta] \E\F\E\To[\heps\;\E] \E$
is the identity,
then the composition of $\F\To[\heta\,\F] \F\E\F\To[\F\,\heps]\F$
is  also the identity.

Thus we have reduced Theorem~\ref{th:main}
to the three statements $\eqref{eq1}$, $\eqref{eq2}$
and $\eqref{eq3}$, which will be proved in the next section.


\section{Proof of the three statements}\label{sec:proof}
\subsection{Intertwiner}

Let us set $\vphi_a\in R(n)$ as follows:
\eqn
\vphi_a e(\nu)
&=&
(x_a\tau_a-\tau_a x_a)e(\nu)=(\tau_ax_{a+1}-x_{a+1}\tau_a)e(\nu)\\
&=&\bl(x_a-x_{a+1})\tau_a+1\br e(\nu)
=\bl\tau_a(x_{a+1}-x_{a})-1\br e(\nu)
\eneqn
if $\nu_a=\nu_{a+1}$ and
$\vphi_a e(\nu)=\tau_ae(\nu)$ if $\nu_a\not=\nu_{a+1}$.
It is called the {\em intertwiner}.

The following lemma is well-known (for example, it easily follows from
by the polynomial representation of \KLRs\ (\cite[Proposition~2.3]{KL09},
\cite[Proposition~3.12]{R08}).
\begin{Lem} \label{lem:ga}
\bnum
\item
For $1 \le a \le n$, we have
$$x_{s_a(b)}\vphi_{a} = \vphi_{a}x_{b} (1\le b\le n+1).$$
\item $\vphi_a^2=Q_{a,a+1}+e_{a, a+1}$.
\item
$\{\vphi_k\}_{1\le k<n}$ satisfies the braid relation.
\item
For $w\in S_n$ and $1\le k<n$,
if $w(k+1)=w(k)+1$, then $\vphi_w\tau_k=\tau_{w(k)}\vphi_w$.
\item
In particular
\eqn
&&\tau_{a} \vphi_{a+1} \vphi_{a} = \vphi_{a+1} \vphi_{a} \tau_{a+1},
\quad\text{and}\quad
\tau_{a+1} \vphi_{a} \vphi_{a+1} = \vphi_{a} \vphi_{a+1} \tau_{a},\\
&&\tau_k\vphi_{a}\cdots\vphi_{n-1}=
\vphi_{a}\cdots\vphi_{n-1}\tau_{k-1}\quad\text{for $a<k\le n-1$.}\eneqn
\enum
\enlemma

\subsection{}
Let  us take $\beta\in Q^+$ with $\haut(\beta)=n$ and $i\in I$.
Let $p$ be the number of times
that $\alpha_i$ appears in $\beta$.
The following lemma is proved by repeated use of 
Theorem~\ref{th:R(n)xR(n)}.

\Lemma\label{lem:dirL}
We have
\eqn
&&e(\beta,i^2)R(\beta+2\alpha_i)e(\beta+\alpha_i,i)
\otimes_{R(\beta+\al_i)}R^\Lambda(\beta+\alpha_i)\\
&&\hs{10ex}
\simeq R(\beta)e(\beta-\al_i,i)\otimes\cor\tau_n\tau_{n+1}
\otimes_{R(\beta-\al_i)}e(\beta-\al_i,i)
R^\Lambda(\beta+\al_i)\\
&&\hs{20ex}\soplus\tau_{n+1}\cor[x_{n+2}]\otimes
e(\beta,i)\RL(\beta+\alpha_i)\\
&&\hs{30ex}\soplus\cor[x_{n+2}]\otimes e(\beta,i)\RL(\beta+\alpha_i).
\eneqn
\enlemma

\Proof
We have
\eqn
\hs{-5ex}&&e(\beta,i^2)R(\beta+2\alpha_i)e(\beta+\alpha_i,i)
\otimes_{R(\beta+\al_i)}R^\Lambda(\beta+\alpha_i)\\
\hs{-5ex}&&
=e(\beta,i^2)\Bigl(R(\beta+\alpha_i)e(\beta,i)\tau_{n+1}
\otimes_{R(\beta)}e(\beta,i)R(\beta+\alpha_i)
\soplus \cor[x_{n+2}]\otimes_\cor R(\beta+\alpha_i)\Bigr)\\
\hs{-5ex}&&\hs{40ex}\otimes_{R(\beta+\al_i)}R^\Lambda(\beta+\alpha_i) \\
\hs{-5ex}&&=e(\beta,i^2)\Bigl(R(\beta)e(\beta-\al_i,i)\tau_{n}
\otimes_{R(\beta-\al_i)}e(\beta-\al_i,i)R(\beta)
\oplus \cor[x_{n+1}]\otimes R(\beta)\Bigr)\tau_{n+1}\\
\hs{-5ex}&&\hs{30ex}\otimes_{R(\beta)}e(\beta,i)\RL(\beta+\alpha_i)\\
\hs{-5ex}&&\hs{15ex}\soplus\cor[x_{n+2}]\otimes_\cor R^\Lambda(\beta+\alpha_i)\\
\hs{-5ex}&&=e(\beta,i^2)R(\beta)e(\beta-\al_i,i)\tau_{n}\tau_{n+1}
\otimes_{R(\beta-\al_i)}e(\beta-\al_i,i^2)R^\Lambda(\beta+\alpha_i)\\
\hs{-5ex}&&\hs{10ex}\soplus\cor[x_{n+1}]\tau_{n+1}\otimes
 e(\beta,i^2)R^\Lambda(\beta+\alpha_i)
\soplus\cor[x_{n+2}]\otimes_\cor e(\beta,i^2) R^\Lambda(\beta+\alpha_i).
\eneqn
Then the lemma follows from
$\cor[x_{n+1}]\tau_{n+1}\oplus\cor[x_{n+2}]
=\tau_{n+1}\cor[x_{n+2}]\oplus\cor[x_{n+2}]$.
\QED

We set 
\eqn
K&\seteq& e(\beta,i^2)R(\beta+2\alpha_i)e(\beta+\alpha_i,i)
\otimes_{R(\beta+\al_i)}R^\Lambda(\beta+\alpha_i)\\
&\simeq&\dfrac{e(\beta,i^2)R(\beta+2\alpha_i)e(\beta+\alpha_i,i)}
{e(\beta,i^2)R(\beta+2\alpha_i)\x R(\beta+\al_i)e(\beta+\alpha_i,i)}.
\eneqn
Then $K$ is an $\bl e(\beta,i^2)R(\beta+2\alpha_i)e(\beta,i^2),
R^\La(\beta+\al_i)\otimes\cor[x_{n+2}]\br$-bimodule.

The preceding lemma says
\eqn
&&K=
R(\beta)\tau_n\tau_{n+1}
e(\beta-\al_i,i^3)R^\Lambda(\beta+\al_i)+
\tau_{n+1}\cor[x_{n+2}]e(\beta,i^2)R^\La(\beta+\alpha_i)\\
&&\hs{35ex}+\cor[x_{n+2}]e(\beta,i^2)R^\La(\beta+\alpha_i).
\eneqn
We define the filtration $\{\Fil_k\}_{k\in\Z}$ of $K$ 
by
\eqn
\Fil_k=\begin{cases}

0&\text{if $k<-1$,}\\
R(\beta)\tau_n\tau_{n+1}e(\beta-\al_i,i^3)
R^\Lambda(\beta+\al_i)+e(\beta,i^2)\tau_{n+1}R^\Lambda(\beta+\alpha_i)
&\text{if $k=-1$,}\\
\Fil_{k-1}+e(\beta,i^2)x_{n+2}^k R^{\Lambda}(\beta+\al_i)
+e(\beta,i^2)\tau_{n+1}x_{n+2}^{k+1} R^\Lambda(\beta+\al_i)
&\text{if $k\ge0$.}
\end{cases}
\eneqn
Note that
$\Fil_k=\Fil_{k-1}+e(\beta,i^2)x_{n+2}^k R^{\Lambda}(\beta+\al_i)
+e(\beta,i^2)x_{n+1}^{k+1}\tau_{n+1} R^\Lambda(\beta+\al_i)$ 
for $k\ge0$.

Recall that $\Gr^\Fil_kK\seteq\Fil_k/\Fil_{k-1}$.
Then we have the following lemma that will be used frequently.

\Lemma \label{lem:Grin}
We have
\bnum
\item
the $\Fil_k$'s are $(R(\beta),R^\La(\beta+\al_i))$-bimodules,
\item $\Fil_kx_{n+2}\subset \Fil_{k+1}$ for any $k$,
\item
the right multiplication of $x_{n+2}$ induces an isomorphism
$\Gr_k^\Fil K\isoto \Gr_{k+1}^\Fil K$ for any $k\ge0$,
\item
$\Ker\bl x_{n+2}\cl \Fil_{-1}\to\Gr^\Fil_0K\br
=R(\beta)\tau_n\tau_{n+1}e(\beta-\al_i,i^3)
R^\Lambda(\beta+\al_i)$.
\ee
\enlemma
\Proof
(i) is obvious.

\noi
(ii) follows from 
\eq&&\tau_n\tau_{n+1}x_{n+2}
=\tau_n(x_{n+1}\tau_{n+1}+1)
=(x_n\tau_n+1)\tau_{n+1}+\tau_n.\label{eq:10}
\eneq

\noi
(iii) follows from Lemma~\ref{lem:dirL}.

\smallskip
\noi
Let us prove (iv).
Set $S\seteq R(\beta)\tau_n\tau_{n+1}e(\beta-\al_i,i^3)
R^\Lambda(\beta+\al_i)$.
Then $Sx_{n+2}\subset \Fil_{-1}+e(\beta,i^2)\RL(\beta+\alpha_i)$
by \eqref{eq:10}.
The homomorphism
$\Fil_{-1}/S\to(\Gr_0^\Fil K)/\bl e(\beta,i^2)\RL(\beta+\alpha_i)\br$
is an isomorphism since
it is isomorphic to $\cor\tau_{n+1}\otimes e(\beta,i)\RL(\beta+\alpha_i)
\isoto[{x_{n+2}}]\cor\tau_{n+1}x_{n+2}\otimes e(\beta,i)\RL(\beta+\alpha_i)$.
\QED
As a corollary of the lemma above, we obtain the following
\Lemma\label{lem:grind}
Let $m\in \Z$ and
let $f(x_{n+2})\in R^\La(\beta+\al_i)\otimes\cor[x_{n+2}]$ be a monic polynomial of degree
$r\ge0$ in $x_{n+2}$ and $u\in K$.
Assume that $uf(x_{n+2})\in \Fil_m$. Then we have
\bnum
\item
if $m\ge r-1$, then $u\in \Fil_{m-r}$,
\item
$u x_{n+2}^{k}\in \Fil_{\max(-1,m-r+k)}$ for any $k\ge0$,
\item
$uf(x_{n+2})\equiv ux_{n+2}^r\mod \Fil_{\max(-1,m-1)}$,
\item
if $m<r-1$, then $u\in   R(\beta)\tau_n\tau_{n+1}e(\beta-\al_i,i^3)
R^\Lambda(\beta+\al_i)$.
\ee
\enlemma
\Proof
(i) It is enough to show that
if $u\in \Fil_k$ and $k>m-r$, then $u\in \Fil_{k-1}$.
For such a $u$ we have $uf(x_{n+2})\in \Fil_m\subset \Fil_{k+r-1}$,
and the injectivity of $\Gr_k^\Fil K\To[{f(x_{n+2})=x_{n+2}^{r}}]\Gr_{r+k}^\Fil K$
implies $u\in \Fil_{k-1}$.

\smallskip
\noi
(ii)\quad 
We have $u x_{n+2}^{k}f(x_{n+2})\in \Fil_{m+k}\subset \Fil_{r+\max(-1,m-r+k)}$.
Hence (i) implies that  $u x_{n+2}^{k}\in \Fil_{\max(-1,m-r+k)}$.

\smallskip
\noi
(iii) follows from (ii).

\smallskip
\noi
(iv) By (ii), $u, ux_{n+2}\in \Gamma_{-1}$. Then the assertion follows from
Lemma~\ref{lem:Grin} (iv).
\QED

Our goal of this subsection is to prove Proposition~\ref{prop:main}
below, and  the following lemma is its starting point.
\Lemma
For $\nu\in I^\beta$ we have, as an element of $K$
\eqn
&&\tau_{n+1}\cdots\tau_1\x\vphi_1\cdots\vphi_{n+1}e(\nu,i^2)
\prod_{\substack{a\le n,\;\nu_a=i}}(x_a-x_{n+2})\\
&&\equiv
-\tau_{n+1}\x[n+2]\prod_{\nu_a\not=i}Q_{i,\nu_a}(x_{n+2},x_a)
e(\nu,i,i)
\mod \Fil_{-1}.
\eneqn
\enlemma

\Proof
We have
$\tau_{n+1}\cdots\tau_1\x\vphi_1\cdots\vphi_{n+1}
=\tau_{n+1}\cdots\tau_1\vphi_1\cdots\vphi_{n+1}\x[n+2]$.
We shall show for $a\le n$
\eq
&&\tau_{n+1}\cdots\tau_a\vphi_a\cdots\vphi_{n+1}\x[n+2]e(\nu,i^2)
\prod_{a\le k\le n,\,\nu_k=i}(x_k-x_{n+2})
\prod_{k<a,\,\nu_k\not=i}Q_{i,\nu_a}(x_{n+2},x_k)\nn\\
&&\hs{10ex}\equiv
\tau_{n+1}\cdots\tau_{a+1}\vphi_{a+1}\cdots\vphi_{n+1}\x[n+2]e(\nu,i^2)
\label{eq:indtauphi}\\
&&\hs{20ex}\prod\limits_{a+1\le k\le n,\,\nu_k=i}(x_k-x_{n+2})
\prod\limits_{k<a+1,\,\nu_k\not=i}Q_{i,\nu_a}(x_{n+2},x_k).\nn
\eneq
If $\nu_a\not=i$, it is obvious.
Assume that $\nu_a=i$.
Then 
\eqn
&&\tau_{n+1}\cdots\tau_{a}\vphi_{a}\cdots\vphi_{n+1}\x[n+2]e(\nu,i^2)
(x_a-x_{n+2})\\
&&=\tau_{n+1}\cdots\tau_{a}(x_{a+1}-x_{a})\vphi_{a}\cdots\vphi_{n+1}\x[n+2]e(\nu,i^2)\\
&&=\tau_{n+1}\cdots\tau_{a+1}(\vphi_a+1)\vphi_{a}\vphi_{a+1}
\cdots\vphi_{n+1}\x[n+2]e(\nu,i^2)\\
&&=\tau_{n+1}\cdots\tau_{a+1}(\vphi_a+1)\vphi_{a+1}\cdots\vphi_{n+1}
\x[n+2]e(\nu,i^2)\\
&&=\tau_{n+1}\cdots\tau_{a+1}\vphi_a\vphi_{a+1}\cdots\vphi_{n+1}\x[n+2]
e(\nu,i^2)\\
&&+\tau_{n+1}\cdots\tau_{a+1}\vphi_{a+1}\cdots\vphi_{n+1}\x[n+2]e(\nu,i^2).
\eneqn
We shall show that
for any $f(x_{n+2})$ and $g=g(x_1,\ldots,x_n)$,
we have
\eq
&&\tau_{n+1}\cdots\tau_{a+1}\vphi_a\vphi_{a+1}\cdots\vphi_{n+1}e(\nu,i^2)
f(x_{n+2})g\in \Fil_{-1}.
\label{eq:vanf}
\eneq
We have
\eqn
&&\tau_{n+1}\cdots\tau_{a+1}\vphi_a\vphi_{a+1}\cdots\vphi_{n+1}e(\nu,i^2)
f(x_{n+2})\\
&&=\tau_{n+1}\cdots\tau_{a+1}f(x_a)\vphi_a\vphi_{a+1}\cdots\vphi_{n+1}e(\nu,i^2)\\
&&=f(x_a)\tau_{n+1}\cdots\tau_{a+1}\vphi_a\vphi_{a+1}\cdots\vphi_{n+1}e(\nu,i^2)\\
&&=f(x_a)\vphi_a\vphi_{a+1}\cdots\vphi_{n+1}\tau_{n}\cdots\tau_{a}e(\nu,i^2).
\eneqn
We have
\eqn
&&\vphi_{n}\vphi_{n+1}=\vphi_n(x_{n+1}\tau_{n+1}-\tau_{n+1}x_{n+1})\\
&&=x_n(x_n\tau_n-\tau_nx_n)\tau_{n+1}-(x_n\tau_n-\tau_nx_n)\tau_{n+1}x_{n+1}\\
&&=x_n^2\tau_n\tau_{n+1}-x_n\tau_n\tau_{n+1}x_n-x_n\tau_n\tau_{n+1}x_{n+1}
+\tau_n\tau_{n+1}x_nx_{n+1}
\eneqn
and it belongs to $\Fil_{-1}$.
Hence we obtain \eqref{eq:vanf}. Then the repeated use of 
\eqref{eq:indtauphi}
implies that
\eqn
&&\tau_{n+1}\cdots\tau_1\vphi_1\cdots\vphi_{n+1}\x[n+2]e(\nu,i^2)
\prod_{k\le n,\,\nu_k=i}(x_k-x_{n+2})\\
&&\equiv
\tau_{n+1}\vphi_{n+1}\x[n+2]e(\nu,i^2)
\prod\limits_{\nu_k\not=i}Q_{i,\nu_a}(x_{n+2},x_k).
\eneqn
Finally 
$\tau_{n+1}\vphi_{n+1}e(\nu,i^2)=\tau_{n+1}(\tau_{n+1}(x_{a+1}-x_a)-1)e(\nu,i^2)
=-\tau_{n+1}e(\nu,i^2)$.

\QED

\Lemma
The following equality holds as an element of $K$.
\eqn
&&\tau_{n+1}\cdots\tau_1\x\vphi_1\cdots\vphi_{n+1}e(\nu,i^2)\\
&&\hs{10ex}=\tau_{n+1}\cdots\tau_1\x\tau_1\cdots\tau_{n+1}e(\nu,i^2)
\hs{-3ex}\prod_{\substack{\text{$k=n+1$ or $\nu_k=i$}}}(x_{n+2}-x_a).\eneqn
\enlemma
\Proof
It is enough to show that
\eq&&\ba{l}
\tau_{n+1}\cdots\tau_1\x
\tau_1\cdots\tau_{a-1}\vphi_{a}\cdots\vphi_{n+1}e(\nu,i^2)\\
\hs{2ex}=\tau_{n+1}\cdots\tau_1\x\tau_1\cdots\tau_{a}\vphi_{a+1}\cdots\vphi_{n+1}
e(\nu,i^2)
(x_{n+2}-x_k)^{\substack{\delta(\text{$a=n+1$ or $\nu_a=i$})}}.
\ea\label{eq:23}
\eneq

If $\nu_a\not=i$ it is trivial. If $\nu_a=i$ or $a=n+1$ then we have 
\eqn
\vphi_{a}\cdots\vphi_{n+1}e(\nu,i^2)
&=&(\tau_a(x_{a+1}-x_a)-1)\vphi_{a+1}\cdots\vphi_{n+1}e(\nu,i^2)\\
&=&\tau_a\vphi_{a+1}\cdots\vphi_{n+1}(x_{n+2}-x_a)e(\nu,i^2)
-\vphi_{a+1}\cdots\vphi_{n+1}e(\nu,i^2).
\eneqn
Since
\eqn&&
\tau_{n+1}\cdots\tau_1\x\tau_1\cdots\tau_{a-1}\vphi_{a+1}\cdots\vphi_{n+1}
e(\nu,i^2)\\
&&\hs{10ex}
=\tau_{n+1}\cdots\tau_1\vphi_{a+1}\cdots\vphi_{n+1}\x\tau_1\cdots\tau_{a-1}
e(\nu,i^2)\eneqn
vanishes as an element of $K$ for $a\le n+1$,
we obtain \eqref{eq:23}.
\QED

Thus we have
\eqn
&&
(-1)^{p}\tau_{n+1}\cdots\tau_1\x\tau_1\cdots\tau_{n+1}e(\nu,i^2)
\prod_{\substack{\text{$a=n+1$ or $\nu_a=i$}}}(x_{n+2}-x_a)^2\\
&&\equiv
-\tau_{n+1}\x[n+2]\prod_{\nu_a\not=i}Q_{i,\nu_a}(x_{n+2},x_a)
e(\nu,i^2)(x_{n+2}-x_{n+1})
\mod \Fil_{-1}.
\eneqn

We have
$\tau_{n+1}(x_{n+2}-x_{n+1})\in \Fil_{0}$,
and hence Lemma~\ref{lem:grind} implies
\eqn
&&\tau_{n+1}\x[n+2]\prod_{\nu_a\not=i}Q_{i,\nu_a}(x_{n+2},x_a)e(\nu,i^2)
(x_{n+2}-x_{n+1})\\
&&\hs{10ex}
\equiv\tau_{n+1}x_{n+2}^{\lan h_i,\La-\beta\ran+2p+1}e(\nu,i^2)
\prod_{\nu_a\not=i}t_{i\nu_a}
\mod \Fil_{\lan h_i,\La-\beta\ran+2p-1}.
\eneqn
In particular 
$$\tau_{n+1}\cdots\tau_1\x\tau_1\cdots\tau_{n+1}e(\nu,i^2)
\prod_{\substack{\text{$a=n+1$ or $\nu_a=i$}}}(x_{n+2}-x_a)^2
\in  \Fil_{\lan h_i,\La-\beta\ran+2p}.$$
Hence
it is equivalent to
$\tau_{n+1}\cdots\tau_1\x\tau_1\cdots\tau_{n+1}e(\nu,i,i)x_{n+2}^{2p+2}$
modulo $\Fil_{\lan h_i,\La-\beta\ran+2p-1}$.

Thus we obtain the following proposition.
\Prop\label{prop:main}
For $\beta\in Q^+$, let $p$ be the number of times that 
$\al_i$ appears in $\beta$, and set
$\la\seteq\La-\beta$, $\la_i\seteq\lan h_i,\la\ran$.
Then there exists $c\in\cor_0*\times$ such that
$$\tau_{n+1}x_{n+2}^{\la_i+2p+1}e(\beta,i^2)
\equiv c\tau_{n+1}\cdots\tau_1\x\tau_1\cdots\tau_{n+1}e(\beta,i^2)x_{n+2}^{2p+2}
\mod \Fil_{\la_i+2p-1}.
$$
\enprop
Note that $\la_i+2p\ge0$.

\subsection{}

Let us define two homomorphisms
\eqn
&&P\cl R(\beta)e(\beta-\al_i,i)\tens_{R(\beta-\al_i)}
e(\beta-\al_i,i^2)\RL(\beta+\al_i)\To K\quad\text{and}\\
&&E\cl R(\beta)e(\beta-\al_i,i)\tens_{R(\beta-\al_i)}
e(\beta-\al_i,i^2)\RL(\beta+\al_i)\To e(\beta,i)\RL(\beta+\al_i)
\eneqn
by $P(a\otimes b)=a\tau_n\tau_{n+1}\otimes b$
and $E(a\otimes b)=ab$.
Then $P$ is injective and Lemma~\ref{lem:grind} implies
\eq
&&\Im(P)=\Ker\bl x_{n+2}\cl \Fil_{-1}\To\Gr_0^\Fil K\br.
\eneq

We can see that $R(\beta)e(\beta-\al_i,i)\tens_{R(\beta-\al_i)}
e(\beta-\al_i,i^2)\RL(\beta+\al_i)$ has a structure of
$(\R(\beta)\otimes\cor\lan x_{n},x_{n+1},\tau_{n}\ran,%
\cor[x_n]\otimes\RL(\beta+\al_i))$-bimodule
by
\eqn
&&(a\otimes b)(x_n\otimes1)=ax_n\otimes b,\\
&&(1\otimes\tau_{n})(a\otimes b)=a\otimes \tau_nb,\\
&&(1\otimes x_k)(a\otimes b)=a\otimes x_kb\quad\text{for $k=n,n+1$.}
\eneqn
Here $\cor\lan x_{n},x_{n+1},\tau_{n}\ran$ is the $\cor$-subalgebra of
$e(\beta-\al_i,i^2)\RL(\beta+\al_i)e(\beta-\al_i,i^2)$
generated by $x_{n},x_{n+1},\tau_{n}$, and it is isomorphic to
the nil affine Hecke algebra $R(2\al_i)$.

\Lemma
For any $z\in R(\beta)e(\beta-\al_i,i)\tens_{R(\beta-\al_i)}
e(\beta-\al_i,i^2)\RL(\beta+\al_i)$, we have
$$P(z)x_{n+2}=P(z(x_n\otimes1))+\tau_{n+1}E(z)+
E((1\tens\tau_{n})z).$$
\enlemma
\Proof
For $z=a\otimes b$, we have
\eqn
P(a\tens b)x_{n+2}&=&a\tau_n\tau_{n+1}x_{n+2}\tens b
=a\tau_n(x_{n+1}\tau_{n+1}+1)\tens b\\
&=&a(x_n\tau_n+1)\tau_{n+1}\tens b+1\tens a\tau_nb.
\eneqn
\QED

\Cor\label{cor:cr}
If $z\in R(\beta)e(\beta-\al_i,i)\tens\limits_{R(\beta-\al_i)}
e(\beta-\al_i,i^2)\RL(\beta+\al_i)$ satisfies
$P(z)x_{n+2}\in\Fil_{-1}$, then $ E\bl(1\tens\tau_{n})z\br=0$.
\encor
Indeed, $P(z)x_{n+2}\equiv E((1\tens\tau_{n})z)\mod \Fil_{-1}$.

\medskip
Set $K^\Lambda=e(\beta, i^2)R^\Lambda(\beta^+2\al_i)e(\beta+\al_i,i)$.
Hence we have
\eqn
K&\simeq& \dfrac{e(\beta, i^2)R(\beta+2\al_i)e(\beta+\al_i,i)}
{e(\beta, i^2)R(\beta^+2\al_i)\x R(\beta+\al_i)e(\beta+\al_i,i)},\\
K^\Lambda&\simeq& \dfrac{e(\beta, i^2)R(\beta+2\al_i)e(\beta+\al_i,i)}
{e(\beta, i^2)R(\beta^+2\al_i)\x R(\beta+2\al_i)e(\beta+\al_i,i)}.
\eneqn
Then there exists a surjective homomorphism
$$p\cl K\to K^\La.$$
Note that
\eq
&& p\bl\tau_{n+1}\cdots\tau_1\x\tau_1\cdots\tau_{n+1}e(\beta,i^2)\br=0.
\eneq

Let us denote by $\{\Fil^\Lambda_k\}_{k\in\Z}$ the filtration of $K^\La$
induced by the filtration $\Fil$ of $K$.

\subsection{Proof of  (\ref{eq1})}
Assume that $\la_i\ge2$.
The statement \eqref{eq1} can be read as
$$
x_{n+2}^{\la_i-2}e(\beta,i^2)+x_{n+1}^{\la_i-1}\tau_{n+1}e(\beta,i^2)
\in \Fil^\La_{\la_i-3}\quad\text{as an element of $K^\La$.}$$
By Proposition~\ref{prop:main} and Lemma~\ref{lem:Grin}, we have
$$
\tau_{n+1}x_{n+2}^{\la_i-1}e(\beta,i^2)
\equiv c\tau_{n+1}\cdots\tau_1\x\tau_1\cdots\tau_{n+1}e(\beta,i^2)
\mod \Fil_{\la_i-3}
$$
as an element of $K$.
Then the desired result holds since
$$\tau_{n+1}x_{n+2}^{\la_i-1}e(\beta,i^2)
\equiv (x_{n+1}^{\la_i-1}\tau_{n+1}+x_{n+2}^{\la_i-2})
e(\beta,i^2)\mod \Fil_{\la_i-3}.$$

\subsection{Proof  of (\ref{eq2})}
Assume that $\la_i=1$.
Set $$w\seteq
\tau_{n+1}e(\beta,i^2)-
c\tau_{n+1}\cdots\tau_1\x\tau_1\cdots\tau_{n+1}e(\beta,i^2)\in K.$$
Then Proposition~\ref{prop:main} together with Lemma~\ref{lem:grind} (iv)
implies that $w,\,wx_{n+2}\in\Fil_{-1}$ and
$w$ belongs to $\Im(P)$.
Hence we can write
$w=P(z)$ for some
$z\in R(\beta)e(\beta-\al_i,i)\tens_{R(\beta-\al_i)}
e(\beta-\al_i,i^2)\RL(\beta+\al_i)$.
Then Corollary~\ref{cor:cr} implies that 
$$E(1\tens\tau_{n})z)=0.$$

Let us define the morphism
$$T
\cl R(\beta)e(\beta-\al_i,i)\kern-1ex\tens_{R(\beta-\al_i)}\kern-1ex
e(\beta-\al_i,i^2)\RL(\beta+\al_i)\to
e(\beta,i)\RL(\beta+\al_i)\tens_{\RL(\beta)}
e(\beta,i)\RL(\beta+\al_i)$$
by $T(a\otimes b)=(a\tau_n)\tens b$.
Then we have
\eqn
\varSigma(T(z))&=&p(P(z)),\\
\tE(T(z))&=&E((1\tens \tau_n)z)=0.
\eneqn

Let us show that $u\seteq e(\beta,i)\tens e(\beta,i)-T(z)$ 
satisfies the condition \eqref{eq2}.

\medskip\noi
(a) \quad $\varSigma(T(z))=p(P(z))=\tau_{n+1}e(\beta,i^2)$
as an element of $\RL(\beta+2\al_i)$.

\smallskip
\noi
(b)\quad $\tE(u)=e(\beta,i)-E\bl T(z)\br
=e(\beta,i)$.

\smallskip
\noi
(c) is obvious.

\subsection{Proof  of (\ref{eq3})}
Assume that $\la_i\le0$. Note that $\ell\seteq-\la_i\le2p$.
Then Proposition~\ref{prop:main} 
says that,
by setting $w\seteq c\tau_{n+1}\cdots\tau_1\x\tau_1\cdots\tau_{n+1}e(\beta,i^2)$,
the element
$\bl wx_{n+2}^{2+\ell}-\tau_{n+1}x_{n+2}e(\beta,i^2)\br x_{n+2}^{-\ell+2p}$
of $K$ belongs to $\Fil_{-\ell+2p-1}$.

Hence we have
$$wx_{n+2}^{\ell+2}-\tau_{n+1}x_{n+2}e(\beta,i^2)
\in \Fil_{-1}.$$
Since $\tau_{n+1}x_{n+2}e(\beta,i^2)\in\Fil_{0}$, we have
$wx_{n+2}^{\ell+2}\in \Fil_{0}$.
Hence Lemma~\ref{lem:grind} implies that
$w x_{n+2}^k\in \Fil_{-1}$ for $0\le k\le \ell+1$.
We set
$$wx_{n+2}^k=P(z_k)+\tau_{n+1}y_k\quad\text{for $0\le k\le \ell+1$}$$
with
$z_k\in R(\beta)e(\beta-\al_i,i)\tens_{R(\beta-\al_i)}
e(\beta-\al_i,i^2)\RL(\beta+\al_i)$ and $y_k\in e(\beta,i)\RL(\beta+\al_i)$.
Then we have for $1\le k\le \ell+2$
\eqn
wx_{n+2}^{k}&=&\bl P(z_{k-1})+\tau_{n+1}y_{k-1}\br x_{n+2}\\
&=&P(z_{k-1}(x_n\tens1))+\tau_{n+1}E(z_{k-1})+E((1\tens\tau_n)z_{k-1})
+\tau_{n+1}x_{n+2}y_{k-1}.
\eneqn
Hence Lemma~\ref{lem:dirL} implies
\eqn
&&z_k=z_{k-1}(x_n\tens1)\quad\text{for $1\le k\le \ell+1$,}\\
&&y_k=E(z_{k-1})\quad\text{for $1\le k\le \ell+1$,}\\
&&E\bl(1\tens\tau_n)z_{k-1}\br=0\quad\text{for $1\le k\le \ell+1$,}\\
&&y_{k-1}=0\quad\text{for $1\le k\le\ell+1$.}
\eneqn
Since $wx_{n+2}^{\ell+2}\equiv \tau_{n+1}x_{n+2}e(\beta,i^2)\mod\Fil_{-1}$,
we have
$y_{\ell+1}=e(\beta,i)$ and $E\bl(1\tens\tau_n)z_{\ell+1}\br=0$.
Thus we obtain $z_k=z_0(x_n^k\tens1)$ for $0\le k\le 1+\ell$,
and 
\eq
&&E\bl z_0(x_n^k\tens1)\br =\begin{cases}
0&\text{for $0\le k\le \ell-1$}\\
e(\beta,i)&\text{for $k=\ell$.}\label{eq:z3}
\end{cases}\\
&&E\bl(1\tens\tau_n)z_{0}(x_n^k\tens1)\br=0
\quad\text{$0\le k\le 1+\ell$.}\label{eq:z4}
\eneq

Let us denote by
\eqn&&q\cl  R(\beta)e(\beta-\al_i,i)\tens_{R(\beta-\al_i)}
e(\beta-\al_i,i^2)\RL(\beta+\al_i)\\
&&\hs{20ex}\To
 \RL(\beta)e(\beta-\al_i,i)\tens_{\RL(\beta-\al_i)}
e(\beta-\al_i,i^2)\RL(\beta+\al_i)
\eneqn
the canonical homomorphism, and
set $v=q(z_0)$.
Then (a) and (b) in \eqref{eq3} follow from
$T_k(v)=E\bl z_0(x_n^k\tens1)\br$.
The equality $G(v)=0$ follows from 
$G(v)=p(P(z_0))=p(w)=0$.

Finally let us prove (d) and (e).
We have
$$E\bl(1\tens x_n^k\tau_n)z_0\br=
E\bl(1\tens\tau_n)z_0(x_n^k\tens1)\br=0\quad\text{for $0\le k\le\ell+1$}$$
by \eqref{eq:z4}.
On the other hand we have
$H_k(v)=E\bl(1\tens \tau_nx_{n+1}^k)z_0\br$.
Since $\tau_nx_{n+1}^k=x_n^k\tau_n+\sum_{a+b=k-1}x_{n+1}^ax_n^b$, we obtain
\eqn
H_k(v)&=&E\bl(1\tens x_n^k\tau_n)z_0\br+\sum_{a+b=k-1}
x_{n+1}^aE\bl(1\tens x_n^b)z_0\br\\
&=&\sum_{a+b=k-1}
x_{n+1}^aE\bl z_0(x_n^b\tens1)\br.
\eneqn
Hence \eqref{eq:z3} implies that
$H_k(v)=0$ for $0\le k\le\ell$ and $H_{\ell+1}(v)=e(\beta,i)$.

\medskip
Thus the proof of \eqref{eq1}, \eqref{eq2} and \eqref{eq1}
is complete.



\bibliographystyle{amsplain}



\end{document}